\documentclass[11pt]{article}
\usepackage{amssymb,amsthm,layout,fancyhdr}

\usepackage{epsfig}

\usepackage{amsmath}

\theoremstyle{definition}
\newtheorem{thm}{Theorem}
\newtheorem{lem}[thm]{Lemma}

\newtheorem{cor}[thm]{Corollary}
\newtheorem{exmp}[thm]{Example}
\newtheorem{defn}[thm]{Definition}

\begin{document}

\title{\bf The
       R\'enyi-Ulam Pathological Liar Game with a Fixed Number of Lies}
\author{%
\sc Robert B.\ Ellis\thanks{Partially supported by NSF grant
    DMS-9977354.} \\
 \small Texas A\&M University,  College Station, Texas\\
 \small{\tt rellis@math.tamu.edu}\\
 \and
\sc Vadim Ponomarenko \\
 \small Trinity University,
 San Antonio, Texas\\
 \small{\tt vadim@trinity.edu}
\and
 \sc Catherine H.\ Yan\thanks{Partially supported by NSF grants
   DMS-0245526 and DMS-0308827 and a Sloan Fellowship. The author is
also affiliated with Dalian University of Technology. }\\
\small Texas A\&M University,
 College Station, Texas\\
\small{\tt cyan@math.tamu.edu}}

\date{Draft of July 19, 2004}

\maketitle

\begin{abstract}
The $q$-round R\'enyi-Ulam pathological liar game with $k$ lies on
the set $[n]:=\{1,\ldots,n\}$ is a 2-player perfect information
zero sum game.  In each round Paul chooses a subset $A\subseteq
[n]$ and Carole either assigns 1 lie to each element of $A$ or to
each element of $[n]\setminus A$. Paul wins if after $q$ rounds
there is at least one element with $k$ or fewer lies.  The game is
dual to the original R\'enyi-Ulam liar game for which the winning
condition is that at most one element has $k$ or fewer lies.  We
prove the existence of a winning strategy for Paul to the
existence of a covering of the discrete hypercube with certain
relaxed Hamming balls.  Defining $F^*_k(q)$ to be the minimum $n$
such that Paul can win the $q$-round pathological liar game with
$k$ lies and initial set $[n]$, we find $F^*_1(q)$ and $F^*_2(q)$
exactly. For fixed $k$ we prove that $F_k^*(q)$ is within an
absolute constant (depending only on $k$) of the sphere bound,
$2^q/\binom{q}{\leq k}$; this is already known to hold
for the original R\'enyi-Ulam liar game due to a result of J.\ Spencer.
\end{abstract}

\section{Introduction}

In this paper we consider the following 2-player perfect information zero-sum
game, which we call the {\em R\'enyi-Ulam pathological liar game}, first
defined in \cite{EY04}. The players Paul and Carole play a $q$-round game on a
set of $n$ elements, $[n]:=\{1,\ldots,n\}$. Each round, Paul splits the set of
elements by choosing a {\em question} set $A\subseteq [n]$; Carole then
completes the round by choosing to assign one {\em lie} either to each of the
elements of $A$, or to each of the elements of $[n]\setminus A$.  A given
element is removed from play, or {\em disqualified}, if it accumulates $k+1$
lies, where $k$ is a predetermined nonnegative constant; in choosing the
question set $A$, we may consider the game to be restricted to the {\em
surviving} elements, which have $\leq k$ lies. The game starts with each
element having no associated lies. If after $q$ rounds at least one element
survives, Paul wins; otherwise Carole wins. Thus Paul plays a strategy to
preserve at least one element for $q$ rounds, and Carole answers adversely. We
think of a capricious or contrary Carole lying ``pathologically'' in order to
disqualify elements as quickly as possible.  Our main result, stated as Theorem
\ref{thm:asymptotic} in Section \ref{sec:vector} and proved in Section
\ref{sec:asymptotic}, is a tight asymptotic characterization of the minimum $n$
for which Paul has a winning strategy for the $q$-round game with a fixed
number, $k$, of lies.

This game arises as the dual to the R\'enyi-Ulam liar game,
originating in \cite{R61} and \cite{U76}, which we refer to as
the {\em original liar game}. The simplest version of the original
game is the ``20 questions'' game in which Paul may ask 20 Yes-No
questions in order to identify a distinguished element $x$ from a
set $[n]$, where Carole answers ``Yes'' or ``No'' without lying.
Here, Paul has a winning strategy iff $\log_2{n}\leq 20$. In the
general version, the number of rounds $q$ and number of elements
$n$ are predetermined, as is the number, $k$, of times Carole is
allowed to lie. We take the equivalent viewpoint that the distinguished
element is not chosen ahead of time by Carole, but rather that she
must answer consistently with there being at least one candidate
for the distinguished element at each round.  Thus a candidate
element $y\in[n]$ cannot be the distinguished element if it would
cause Carole to have lied about it $k+1$ times.  Paul's strategy
in the original game, therefore, is to win by forcing Carole to
associate $k+1$ lies  with all but one element within $q$ rounds,
and Carole's strategy is to answer questions adversely so that at
least two candidate elements remain after $q$ rounds. Recently,
Pelc thoroughly surveyed what is known about the original liar
game and many of its variants \cite{P02}.

The duality between the pathological liar game and the original
liar game arises from the choice of Paul's condition to win.  In
the pathological liar game at least one element must survive for
Paul to win, but in the original game at most one element may
survive for him to win. The remaining mechanics of the two games
are the same, in that each round Paul chooses a question subset
$A\subseteq [n]$ and Carole decides to assign lies either to $A$
or to $[n]\setminus A$.

In Section \ref{sec:vector}, we describe how each stage of the
pathological game can be encoded in a $(k+1)$-tuple state vector
which keeps track of the number of lies associated with each
element. In Section \ref{sec:Berlekamp} we discuss the Berlekamp
weight function on a state vector and how  a winning strategy by
Paul corresponds to maximizing (minimizing) the weight of the
state vector after $q$ rounds in the pathological (original) liar
game. In Section \ref{sec:asymptotic}, we give the value of $n$,
up to a constant independent of $q$, for which Paul can win the
$q$-round game with a fixed number, $k$, of lies. In Sections
\ref{sec:oneLie} and \ref{sec:twoLie}, we give the exact minimum
$n$ for which Paul can win the $q$-round 1-lie and 2-lie games,
respectively. Finally, in Section \ref{sec:equiv}, we prove the
equivalence of the existence of a winning strategy for Paul in the
pathological (original) liar game to the existence of a covering
(packing) of the hypercube with certain relaxed Hamming spheres,
and discuss the connection to covering codes and error-correcting
codes.

\section{The vector game format\label{sec:vector}}

The mechanics of both the pathological liar game and the original
liar game are encapsulated in the following vector framework due
to Berlekamp \cite{B68}. Given that the game parameters are $n$
elements, $q$ rounds, and $k$ lies, the initial state of the game
is the $(k+1)$-vector $(n,0,\ldots,0)$.  An intermediate stage of
the game after some number of rounds is encoded by the {\em state
vector} $\vec{x}=(x_0,x_1,\ldots,x_k)$, where $x_i$ denotes the
number of elements of $[n]$ associated with $i$ lies (disqualified
elements, with $k+1$ lies, are not tracked by the state vector).
The state vector completely encodes a stage of the game because an
element of $[n]$ is distinguished only by the number of lies
associated with it. Paul chooses a question set $A\subseteq [n]$
corresponding to an integer {\em question vector}
$\vec{a}=(a_0,a_1,\ldots,a_k)$ which must be {\em legal}, that is,
$0\leq a_i\leq x_i$ for each $i\in\{0,\ldots,k\}$.
Carole answers either ``Yes'' or ``No.'' By
answering ``Yes,'' Carole assigns an additional lie to each
element in $[n]\setminus A$, so that the next state vector
$Y(\vec{x},\vec{a})$ is obtained from $\vec{x}$ by moving elements
corresponding to $[n]\setminus A$ to the right one position.
Analogously, by answering ``No,'' Carole causes the next state
vector $N(\vec{x},\vec{a})$ to arise from moving elements
corresponding to $A$ to the right one position. Therefore the
subsequent state chosen by Carole is either
\begin{equation}
\begin{array}{rcrrrcrl}
Y(\vec{x},\vec{a}) & := & ( & a_0, & a_1+x_0-a_0, & \ldots, &
    a_k+x_{k-1}-a_{k-1}) & \quad\mbox{or}\\
N(\vec{x},\vec{a}) & := & ( & x_0-a_0, & x_1-a_1+a_0,
    & \ldots, & x_k-a_k+a_{k-1}). &
\end{array} \label{eqn:nextState}
\end{equation}
Elements which become associated with $k+1$ lies are considered to
be shifted out of the state vector to the right, and so we may
consider the question set $A$ and the set of elements $[n]$ to be
restricted at any given stage to the surviving elements. In the
pathological liar game, Paul wins iff after $q$ rounds
$\sum_{i=0}^k x_i\geq 1$ (at least one element survives).
In the original liar game, Paul wins iff
after $q$ rounds $\sum_{i=0}^k x_i\leq 1$.

More generally, we may consider a game starting with an arbitrary
nonnegative state vector $\vec{x}=(x_0,\ldots,x_k)$.  We will use the
following shorthand.
\begin{defn}
(i) The {\em $(\vec{x},q,k)^*$-game} is the $q$-round pathological
liar game with $k$ lies and initial state $\vec{x}$.\\
(ii) The {\em $(\vec{x},q,k)$-game} is the $q$-round original liar game
with $k$ lies and initial state $\vec{x}$.\\
In either game, the initial state $\vec{x}=(x_0,\ldots,x_k)$ encodes
for $0\leq i\leq k$ the number $x_i$ of elements which are initially
associated with $i$ lies.

\end{defn}
The $k$ is redundant when $\vec{x}$ is specified. Both games are
monotonic in the following sense.  Suppose $\vec{x}=(x_0,\ldots,x_k)$,
$\vec{y}=(y_0,\ldots,y_k)$, and $0\leq y_i\leq x_i$ for all $0\leq i\leq k$;
i.e., $\vec{x}$ {\em covers} $\vec{y}$.
If Paul has a strategy to win the
$(\vec{y},q,k)^*$-game (the $(\vec{x},q,k)$-game),
then he has a strategy to win the $(\vec{x},q,k)^*$-game
(the $(\vec{y},q,k)$-game). The new strategy is obtained
from the winning strategy in the pathological game by arbitrarily
choosing whether the extra elements corresponding to $x_i-y_i$
are in $A$ or $[n]\setminus A$,
and in the original game by restricting all questions $A$ by
intersection with the set of all elements represented by
$y_0,\ldots,y_k$. In fact, the same monotonicity holds if
$\vec{x}$ {\em majorizes} $\vec{y}$; i.e., if for all $0\leq j\leq k$,
$\sum_{i=0}^j y_i \leq \sum_{i=0}^j x_i$.  Empirically, an element
lasts longer in the game if it starts with fewer associated lies.
Monotonicity under majorization is an immediate result of Theorem
\ref{thm:equivCovering}, as we will describe in Section
\ref{sec:equiv}. We may now define $F^*_k(q)$ to
be the minimum number $n$ such that Paul has a winning strategy
for the $((n,0,\ldots,0),q,k)^*$-game. The previously defined
maximum $n$ such that Paul can win the $((n,0,\ldots,0),q,k)$-game
is $F_k(q)$. Pelc determined $F_1(q)$ exactly in \cite{P87},
Guzicki determined $F_2(q)$ in \cite{G90}, Deppe determined
$F_3(q)$ in \cite{D00}, and Spencer determined $F_k(q)$ for fixed
$k$ to within a constant independent of $q$.
Of particular importance to this paper is the following result of
Spencer, given implicitly in Section 3 of \cite{S92}.

\begin{thm}[Spencer]\label{thm:spencerMain}
For any fixed nonnegative integer $k$ there exist constants
$q_k, C_k$ such that for
all $q\geq q_k$,
$$ \frac{2^q}{\binom{q}{\leq k}}-C_k \ \leq \
    F_k(q) \ \leq \
    \frac{2^q}{\binom{q}{\leq k}}.$$
\end{thm}

Here, $\binom{q}{\leq k}:=\sum_{i=0}^k\binom{q}{i}$ is the size
of a radius $k$ Hamming ball in the $q$-dimensional discrete
hypercube $Q_q$ (Section \ref{sec:equiv} explores this further).
The main result of this paper, which we prove in Section
\ref{sec:asymptotic}, is the following dual of Theorem
\ref{thm:spencerMain}.

\begin{thm}\label{thm:asymptotic}
For any fixed nonnegative integer $k$ there exist constants $q^*_k,
C^*_k$ such that for all $q\geq q^*_k$,
$$ \frac{2^q}{\binom{q}{\leq k}} \ \leq \
    F^*_k(q) \ \leq \
    \frac{2^q}{\binom{q}{\leq k}}+C^*_k.$$
\end{thm}

\section{The Berlekamp weight function\label{sec:Berlekamp}}

For a nonnegative integer $q$ and a state vector
$\vec{x}=(x_0,\ldots,x_k)$, the {\em $q$-weight} of $\vec{x}$ is defined
to be
\begin{equation}\label{eqn:qWeight}
  wt_q(\vec{x}) := \sum_{i=0}^k x_i \binom{q}{\leq k-i}.
\end{equation}
This is the {\em Berlekamp weight function} introduced in
\cite{B68}. The number of ways to select positions
for at most $k-i$ lies in a sequence of Y/N responses by Carole
of length $q$ is
$\binom{q}{\leq k-i}$, which motivates the weight of an element
counted by $x_i$. We will abuse notation and denote
$wt_q((x_0,\ldots,x_k))$ by $wt_q(x_0,\ldots,x_k)$.
We will
see that Carole can always win the $(\vec{x},q,k)^*$-game when
$wt_q(\vec{x})<2^q$. Intuitively, elements with fewer associated
lies are worth more toward a win by Paul. To borrow an analogy
from \cite{S92}, we can think of the $x_i$'s as representing {\em
coins} of various denominations, where we call the coins with
smallest weight, counted by $x_k$, {\em pennies}. We now present a
well-known conservation lemma concerning the weight function,
previously appearing in \cite{B68}.

\begin{lem}[Conservation of weight]\label{lem:conserve}
Let $q\geq 1$, let $\vec{x}$ be a state vector, and let $\vec{a}$ be a legal
question for $\vec{x}$.
Then
$$ wt_q(\vec{x}) \ = \ wt_{q-1}(Y(\vec{x},\vec{a}))
    +wt_{q-1}(N(\vec{x},\vec{a})).  $$
\end{lem}
\begin{proof}
Using (\ref{eqn:nextState}) and (\ref{eqn:qWeight}), we compute
\begin{align}
wt_{q-1} & (Y(\vec{x},\vec{a})) +wt_{q-1}(N(\vec{x},\vec{a}))
    \ = \ x_0\binom{q-1}{\leq k}
    +\sum_{i=1}^k (x_i+x_{i-1})\binom{q-1}{\leq k-i}
    \nonumber\\
& = \sum_{i=0}^k x_i\left(\binom{q-1}{\leq k-i}+
  \binom{q-1}{\leq k-i-1}\right) \ = \ wt_q(\vec{x}), \nonumber
\end{align}
by repeated use of the identity $\binom{n}{k}=\binom{n-1}{k}+
\binom{n-1}{k-1}$.
\end{proof}

The lemma illustrates that Carole's choice in answering ``Yes'' or
``No'' to a question by Paul induces a choice of weight of the
resulting state vector.  In particular, Carole might always choose
the resulting state with lower weight, giving a
constraint on Paul's ability to win the $(\vec{x},q,k)^*$-game
which holds for any $k$. We call the following lemma the
{\em sphere bound} because of a connection to the sphere bound
of coding theory to be made clear after Theorem \ref{thm:equivCovering}.

\begin{lem}[Sphere bound]\label{lem:lowerBound}
Let $q,k\geq 0$ and let $\vec{x}=(x_0,\ldots,x_k)$ be a nonnegative
vector. If $wt_q(\vec{x})<2^q$, then Carole can win the
$q$-round pathological liar game with $k$ lies and initial state
$\vec{x}$. Consequently, $F^*_k(q)\geq 2^q/\binom{q}{\leq k}$.
\end{lem}
\begin{proof}
Regardless of Paul's initial question, by Lemma \ref{lem:conserve}
Carole may respond so that the resulting state has weight at most
$wt_q(\vec{x})/2<2^{q-1}$.  By induction, Carole may respond to
Paul's remaining $q-1$ questions to ensure the 0-weight of the
final state is $<1$. Since the state vector must always be
integer, Carole can always force the vector $(0,\ldots,0)$ in $q$
rounds.
\end{proof}

In the original game, the analog to the above lemma is that Carole
has a strategy to win the $(\vec{x},q,k)$-game when
$wt_q(\vec{x})>2^q$.  This is proved in \cite{S92} by showing that
if Carole answers randomly at each stage, the probability that the
final weight is $>1$  is nonzero, and thus Carole has a winning
strategy since it is a perfect information game. The proof of
Lemma \ref{lem:lowerBound} could be rewritten from this randomized
perspective.

Lemma \ref{lem:lowerBound} shows that a necessary condition for
Paul to win the $q$-round pathological liar game with starting
state $\vec{x}$ is that $wt_q(\vec{x})\geq 2^q$, but in general
this is not sufficient.  Paul is not always able to choose a
question which balances the weights of the possible next states.
Given some intermediate state $\vec{x}$ with $j+1$ rounds
remaining and a question $\vec{a}$, the resulting {\em weight
imbalance} between possible next states is defined
as (cf.\ Section 2 of \cite{S92})
\begin{equation}\label{eqn:imbalance}
\Delta_j(\vec{x},\vec{a}) \ := \ wt_j(Y(\vec{x},\vec{a}))-
   wt_j(N(\vec{x},\vec{a})).
\end{equation}
The following is a counterexample to the converse of Lemma
\ref{lem:lowerBound}.

\begin{exmp}\label{eg:imbalance}
Let $\vec{x}=(3,1)$ be the initial state of a
$((3,1),4,1)^*$-game.  Note that $wt_4((3,1))=3\cdot 5+1\cdot
1=16$, and so Paul could possibly have a winning strategy.  But
any first-round question $\vec{a}$ by Paul will satisfy
$|\Delta_3(\vec{x},\vec{a})|\geq 2$.  One question minimizing
$|\Delta_3(\vec{x},\vec{a})|$ is $\vec{a}=(1,1)$, for which
$Y(\vec{x},\vec{a})=(1,3)$, $N(\vec{x},\vec{a})=(2,1)$, and
$\Delta_3(\vec{x},\vec{a})=7-9=-2$. In any event, Carole responds
so that the next state has 3-weight at most 7, guaranteeing
herself to win the game.
\end{exmp}

Paul's goal in the pathological liar game, in terms of the weight
function, corresponds to maximizing the 0-weight of the game state
after $q$ rounds. The capability to identify
situations in which he can choose ``perfectly balancing'' questions at
every stage so that  $\Delta_j(\vec{x},\vec{a})=0$ would provide
a partial converse to Lemma \ref{lem:lowerBound}; however, this
is sometimes impossible (cf.\ Example \ref{eg:imbalance}), and
difficult to know if it is possible when initially
the $q$-weight is close to $2^q$.

\section{Asymptotics of the $k$-lie game\label{sec:asymptotic}}

Since the full converse to Lemma \ref{lem:lowerBound} is impossible,
we instead wish to identify the states $\vec{x}$ having $wt_q(\vec{x})$
close to $2^q$ for which Paul can win the $(\vec{x},q,k)^*$-game.
As Spencer proved in \cite{S92}, there
is a  large category of states $\vec{x}=(x_0,\ldots,x_k)$ such that if
$wt_q(\vec{x})=2^q$ and $x_k$ is  large enough, then Paul can find $q$
questions which make the weight imbalance vanish {\em at each stage}.
Intuitively two processes are at work. If there are enough
``pennies,'' counted by $x_k$, then $\vec{a}$ can be chosen so that
the weights of the two possible next states $\mbox{Y}(\vec{x},\vec{a})$
and $\mbox{N}(\vec{x},\vec{a})$ are exactly equal.  The number of
pennies in the next state is maintained sufficiently by drawing from
$x_{k-1}$ and $x_k$. To employ Spencer's result, it will
suffice to begin with $\vec{x}$ having $q$-weight slightly more than
$2^q$ and reduce in $k$ rounds to a state $\vec{y}$ with
$(q-k)$-weight exactly $2^{q-k}$ for which Spencer's theorem holds.
Here now is Spencer's result, essentially
appearing as the ``Main Theorem'' in Section 2 of \cite{S92}, in a
form convenient for our purposes.

\begin{thm}[Spencer] Let $k$ be fixed.
There are constants $c, q_0$ (dependent on $k$) so that the
following holds for all $q\geq q_0$: if $wt_q(x_0,\ldots,x_k)=2^q$
and $x_k>cq^k$, then Paul has a strategy to reach a state
$\vec{z}$ with $wt_0(\vec{z})=1$ in exactly $q$ rounds such that
every intermediate state $(u_0,\ldots,u_k)$ after playing $j$
rounds satisfies
$wt_{q-j}(u_0,\ldots,u_k)=2^{q-j}$.\label{thm:spencer}
\end{thm}

\begin{thm}\label{thm:spencerExtend}
Let $k$ be fixed. There are constants $c_1, q_k^*$ (dependent on
$k$) so that the following holds for all $q\geq q_k^*$: if
$wt_q(x_0,\ldots,x_k)\geq 2^q+c_1 \binom{q}{k}$, then Paul can win
the $q$-round pathological liar game with $k$ lies and initial
state $\vec{x}=(x_0,\ldots,x_k)$.
\end{thm}
\begin{proof}
The proof proceeds in three main stages.  First, the first $k$ rounds
of the game are played with a ``floor-ceiling'' question strategy
which ensures that the resulting state $\vec{y'}$ satisfies
$wt_{q-k}(\vec{y'})\geq 2^{q-k}$.  Second, coins are removed from
$\vec{y'}$ to obtain $\vec{y}$ with $(q-k)$-weight exactly $2^{q-k}$.
Finally, Theorem \ref{thm:spencer} is applied to $\vec{y}$ to reach
a state $\vec{z}$ with $wt_0(\vec{z})=1$ after an additional
$q-k$ rounds.

Paul plays the first $k$ rounds of the game, reaching the state
$\vec{y'}=(y'_0,\ldots,y'_k)$,  according to the following
strategy which is oblivious to Carole's responses.  If
$\vec{u}(j)=(u_0(j),\ldots,u_k(j))$ is the state when $j$ rounds
remain, then for $q \geq j > q-k$, Paul's next question
$\vec{a}(j)=(a_0(j),\ldots,a_k(j))$ is defined by letting
$a_i(j)=\lfloor u_i(j)/2 \rfloor$ or $\lceil u_i(j)/2\rceil$, so
that the least $i$ for which $u_i(j)$ is odd results in choosing
$a_i(j)=\lceil u_i(j)/2 \rceil$, and the overall choice of floors and
ceilings for the odd $u_i(j)$'s alternates.

By combining (\ref{eqn:nextState}) and (\ref{eqn:qWeight}) with
the definition of $\Delta_j$ in (\ref{eqn:imbalance}), the weight
imbalance of the two possible next states when $j+1$ rounds remain
is at most
\begin{equation}\label{eqn:nonnegImbalance}
\Delta_{j}(\vec{u}(j+1),\vec{a}(j+1))  =
   \sum_{i=0}^k (2a_i(j+1)-u_i(j+1))\binom{j}{k-i} \leq
   \binom{j}{k} 
   ,
\end{equation}
where we know the value is nonnegative by definition of $\vec{a}(j+1)$.
By Lemma \ref{lem:conserve} and (\ref{eqn:nonnegImbalance}), we have
for each intermediate state $\vec{u}(j+1)$ (with indexes $j+1$
suppressed for clarity)
$$ wt_j(\mbox{Y}(\vec{u},\vec{a})) \geq
   wt_j(\mbox{N}(\vec{u},\vec{a}))
   \geq \frac{wt_{j+1}(\vec{u})-\binom{j}{k}}{2}.
$$
Therefore with an initial state of weight
\begin{equation}
wt_q(\vec{x}) \geq 2^q+c_1 \binom{q}{k}
    \geq 2^q + \sum_{j=q-1}^{q-k} 2^{q-1-j} \binom{j}{k}, \nonumber
\end{equation}
for some constant $c_1$ and $q\geq q_1$ large enough, Paul can
guarantee a state $\vec{y'}$ with $wt_{q-k}(\vec{y'})\geq 2^{q-k}$
after $k$ rounds.

The number of pennies $y'_k$ after $k$ rounds is large, by the
following  argument.  Since $wt_q(\vec{x})\geq 2^q$ and
the largest weight of an element is $\binom{q}{\leq k}\leq q^k$,
then $\sum_{i=0}^k x_i \geq 2^q/q^k$. Thus there exists a coordinate
$i_0$ for which $x_{i_0}\geq 2^q/\big((k+1)q^k\big)$.  By definition
of the first $k$ questions,
\begin{eqnarray}
y'_k & = & u_k(q-k) \geq \lfloor 2^{-1}u_k(q-k+1)\rfloor
    \geq \cdots \geq \lfloor 2^{-i_0} u_k(q-k+i_0)\rfloor \nonumber \\
& \geq & \lfloor 2^{-i_0-1} u_{k-1}(q-k+i_0+1)\rfloor\geq \cdots \geq
    \lfloor 2^{-k} u_{i_0}(q) \rfloor =x_{i_0}\nonumber \\
& \geq &
\left\lfloor 2^{-k}\cdot\frac{2^q}{(k+1)q^k}\right\rfloor
    \geq c_2q^k. \nonumber
\end{eqnarray}
The first line is true because $u_k(j)$ is at least $\lfloor
u_k(j+1)/2\rfloor$, the second line is true because $u_i(j)$ is at
least $\lfloor u_{i-1}(j+1)/2\rfloor$, and the last inequality is
true for any choice of $c_2$ and $q\geq q_2$ provided $q_2$ is
taken to be large enough. We note that the choice of $c_1$ does
not affect the choice of $c_2$ in this analysis.

Now obtain the state $\vec{y}=(y_0,\ldots,y_k)$
with $(q-k)$-weight $2^{q-k}$ from
$\vec{y'}$ by greedily removing coins of decreasing weight, so
that either only $2^{q-k}$ pennies are left, or fewer than
$\binom{q-k}{\leq k}$ pennies were removed. In the first case Paul
trivially can make the game last another $q-k$ rounds; in the
second case at least
$$
y_k \geq c_2q^k-\binom{q-k}{\leq k} \ \geq c_3 (q-k)^k
$$
pennies remain. The constant $c_3$ can be chosen to be at least
$c_2-1$, for instance, provided that $q\geq q_3$ for $q_3$ large
enough. Choose $c_3$ and $q_k^*\geq \max\{q_1,q_2,q_3\}$ large
enough so that $c_3$ and $q_k^*-k$ satisfy the requirements of
Theorem \ref{thm:spencer} for the $(\vec{y},q-k,k)$-game.
Therefore Paul can win the $(\vec{x},q,k)^*$-game.
\end{proof}

\begin{proof}[Proof of Theorem \ref{thm:asymptotic}]
From Lemma \ref{lem:lowerBound},
$F^*_q(k)\geq 2^q/\binom{q}{\leq k}$. Now suppose $q\geq q_k^*$
and let
$n=\lceil (2^q + c_1\binom{q}{k})/\binom{q}{\leq k}\rceil$, where
$c_1$ and $q_k^*$ are  as in Theorem \ref{thm:spencerExtend}. Then
$wt_q(n,0,\ldots,0)\geq 2^q + c_1\binom{q}{k}$ and $F^*_k(q)\leq n
\leq \lceil (2^q+c_1\binom{q}{k})/\binom{q}{\leq k}\rceil \leq
2^q/\binom{q}{\leq k} +C^*_k$ for $q\geq
q^*_k$ and some constant $C^*_k$.
\end{proof}

We remark that the excess weight above $2^q$ in Theorem 8 is
needed so that Paul can guarantee a $(q-k)$-weight of $2^{q-k}$
after the first $k$ rounds and go on to win when $q$ is large
enough.  The exact excess required is difficult to compute
for general $k$.  However, in the next two sections we
will compute the exact amount required for $k=1$ and $2$ for
any $q$, not just when $q$ is large enough.

\section{Exact result for the 1-lie game\label{sec:oneLie}}

We now consider the $q$-round
pathological liar game with 1 lie and initial state $(n,0)$.
For this section, define
the {\em character} $ch(x_0,x_1)$ of a state $(x_0,x_1)$ to be the
maximum $q$ such that $wt_{q}(x_0,x_1)\geq 2^q$.  Furthermore,
denote by $(y_0,y_1)$ the game state immediately following the
state $(x_0,x_1)$ and Paul's question $(a_0,a_1)$, so that
$(y_0,y_1)=(a_0,a_1+x_0-a_0)$ or $(x_0-a_0,a_0+x_1-a_1)$,
depending on Carole's response of ``Y'' or ``N,'' respectively.
The next theorem completely characterizes the values of $n$ for
which Paul can win the $((n,0),q,1)^*$-game.

\begin{thm}\label{thm:1Lie}Let $q\geq 0$.  Paul has a winning strategy
for the $q$-round pathological liar game with 1 lie and initial
state $(n,0)$ iff
\begin{equation} 2^q \leq \left\{\begin{array}{rl}
      n(q+1) & \mbox{if $n$ is even,}\\
      n(q+1)-(q-1) & \mbox{if $n$ is odd.}
    \end{array}\right.
\label{eqn:1LieCondition}
\end{equation}
\end{thm}

The difference in the even and odd cases reflects the fact that when
$n$ is odd, Paul's first question is forced to be inefficient, as
there is no way to balance $a_0$ with $x_0-a_0$.  By considering
the possibilities for $\lceil 2^q/(q+1)\rceil\mbox{ mod }2$ and
$2^q\mbox{ mod }q+1$, it is not difficult to obtain the
following.
\begin{cor}
Let $SB^*_1:=\lceil 2^q/(q+1)\rceil$ be the sphere bound for
the $((n,0),q,1)^*$-game.  Then
$$
F^*_1(q) = \left\{\begin{array}{rl}
    SB^*_1, & \mbox{if $SB^*_1$ is odd and
      $(2^q\mbox{ mod }q+1)\in\{1,2\}$,} \\
    2\lceil SB^*_1/2\rceil, & \mbox{otherwise.}
\end{array}\right.
$$
\end{cor}
The proof of Theorem \ref{thm:1Lie} follows in one direction
by Lemma \ref{lem:1LieLB}, and the other direction will be proved
after Lemmas \ref{lem:smallMSB} and \ref{lem:manyPennies}.  This
proof technique is based on that of Pelc's theorem in Section 2 of
\cite{P87}, which states that the characterization for Paul having a
winning strategy for the $((n,0),q,1)$-game is obtained from
(\ref{eqn:1LieCondition}) by reversing the inequality.

\begin{lem}\label{lem:1LieLB}
Let $q\geq 0$.  Carole can win the $q$-round pathological liar
game with 1 lie and initial state $(n,0)$ provided
$$ 2^q > \left\{\begin{array}{rl}
n(q+1) & \mbox{if $n$ is even,}\\
n(q+1)-(q-1) & \mbox{if $n$ is odd.}
\end{array}\right.
$$
\end{lem}
\begin{proof}
The case of $n$ even follows directly from Lemma
\ref{lem:lowerBound}, since $wt_q(n,0)=n(q+1)$.
If $n$ is odd, observe that whatever Paul's first
question is, Carole may respond so that in the resulting state
$(y_0,y_1)$, $y_0<y_1$, and so
$$ wt_{q-1}(y_0,y_1)\leq \frac{n-1}{2}q+\frac{n+1}{2}=
    \frac{n(q+1)-(q-1)}{2} < 2^{q-1}. $$
Now apply Lemma \ref{lem:lowerBound} to show that Carole can win
the $((y_0,y_1),q-1,1)^*$-game.
\end{proof}

The next lemma handles the late rounds of the game for which there
is at most 1 element with no accumulated lies.

\begin{lem}\label{lem:smallMSB}
Paul can win the $q$-round pathological liar game with 1 lie and
initial state $(x_0,x_1)$ provided $0\leq x_0\leq 1$ and $q\leq
ch(x_0,x_1)$.
\end{lem}
\begin{proof}
Without loss of generality, assume $q=ch(x_0,x_1)$. We prove the
lemma by induction on $q$, by exhibiting a question Paul can ask
that will not reduce the character by more than one. Since
$q=ch(x_0,x_1)$, $wt_q(x_0,x_1)=(q+1)x_0+x_1\ge 2^q$.

If $x_0=0$, then $x_1=wt_q(x_0,x_1)\geq 2^q$; if Paul chooses the
question $\vec{a}=(0,\lfloor\frac{x_1}{2}\rfloor)$, then
$y_1=\lfloor\frac{x_1}{2}\rfloor$ or $\lceil\frac{x_1}{2}\rceil$.
In either case, $wt_{q-1}(y_0,y_1)\geq \lfloor
\frac{2^q}{2}\rfloor \geq 2^{q-1}$, and so $ch(y_0,y_1)\geq q-1$.

If $x_0=1$, set $a_0=1$ and $a_1=\lfloor \frac{x_1+1-q}{2}
\rfloor$. Observe that $a_1\ge 0$, since otherwise $q>x_1+1$ and
$2q>q+1+x_1=wt_q(1,x_1)\ge 2^q$, which is impossible.  Paul then
asks $\vec{a}=(1,a_1)$, and Carole can choose between $(1,a_1)$ or
$(0,x_1+1-a_1)$.  The weight imbalance is
$|\Delta_{q-1}(\vec{x},\vec{a})|=
|(q+a_1)-(x_1+1-a_1)|=|q-x_1-1+2a_1|\le 1$. By Lemma
\ref{lem:conserve} and because $2^q$ is even, we have
$wt_{q-1}(y_0,y_1)\geq 2^{q-1}$. Hence $ch(y_0,y_1)\ge
q-1$.
\end{proof}

We now show that certain state vectors $(x_0,x_1)$ in the game
allow Paul a question which guarantees that the next state has
three narrow constraints, including a character reduced by at most
one.

\begin{lem}\label{lem:manyPennies}
Let $(x_0,x_1)$ be a state with $ch(x_0,x_1)\geq 1$ and
$x_1\ge x_0-1\ge 1$.  Then there exists a question $(a_0,a_1)$
such that regardless of Carole's answer the next state $(y_0,y_1)$
will satisfy:
\begin{align}
&\lfloor \frac{x_0}{2} \rfloor   \le  y_0  \le  \lceil
\frac{x_0}{2} \rceil \label{eqn:Lem1}\\
&y_1  \ge  y_0-1 \label{eqn:Lem2}\\
&ch(y_0,y_1)  \ge  ch(x_0,x_1)-1. \label{eqn:Lem3}
\end{align}
\end{lem}
\begin{proof}
Without loss of generality, assume $q=ch(x_0,x_1)$.
The proof depends on whether $x_0$ is even or odd. {\em Case 1}
($x_0$ is even). Paul chooses the legal question
$\vec{a}=(\frac{x_0}{2},\lfloor\frac{x_1}{2}\rfloor)$ so that
$(y_0,y_1)=(\frac{x_0}{2},\frac{x_0}{2}+\lfloor\frac{x_1}{2}\rfloor)$
or $(\frac{x_0}{2},\frac{x_0}{2}+\lceil\frac{x_1}{2}\rceil)$.
Regardless of Carole's response, $y_0=\frac{x_0}{2}$, satisfying
condition (\ref{eqn:Lem1}); also, $y_1\geq
\frac{x_0}{2}+\lfloor\frac{x_1}{2}\rfloor\geq y_0-1$, satisfying
condition (\ref{eqn:Lem2}). Finally, since $2^q$ is even and
$|\Delta_{q-1}(\vec{x},\vec{a})|
=\lceil\frac{x_1}{2}\rceil-\lfloor\frac{x_1}{2}\rfloor\leq 1$, we
have $wt_{q-1}(y_0,y_1)\geq 2^{q-1}$, and so condition
(\ref{eqn:Lem3}) is satisfied.

{\em Case 2} ($x_0$ is odd). Paul chooses
$\vec{a}=(\frac{x_0+1}{2},\lceil\frac{x_1-q+1}{2}\rceil)$, so that
$(y_0,y_1)=(\frac{x_0+1}{2},\frac{x_0-1}{2}+\lceil
\frac{x_1-q+1}{2}\rceil)$ or
$(\frac{x_0-1}{2},\frac{x_0+1}{2}+x_1-\lceil\frac{x_1-q+1}{2}\rceil)$.
To show the question is legal, we require
$a_1=\lceil\frac{x_1-q+1}{2}\rceil\geq 0$.  Otherwise,
$x_1-q+1<-1$, or $x_1\leq q-3$, and so $x_0\leq q-2$.
With this assumption on $x_0$ and $x_1$, $2^q\leq
wt_q(x_0,x_1)\leq (q+1)(q-2)+q-3=q^2-5$, which is impossible for
$q\geq 0$, and so the question is legal. Continuing, clearly
condition (\ref{eqn:Lem1}) holds.  If Carole answers ``Y,''
$y_1-y_0+1=\lceil\frac{x_1-q+1}{2}\rceil$, which is at least 0. If
Carole answers ``N,''
$y_1-y_0+1=x_1-\lceil\frac{x_1-q+1}{2}\rceil+2$, which is clearly
nonnegative.  Thus condition (\ref{eqn:Lem2}) holds. Again, $2^q$
is even, and $|\Delta_{q-1}(\vec{x},\vec{a})|
=|2\lceil\frac{x_1-q+1}{2}\rceil-(x_1-q+1)|\leq 1$; therefore
$wt_{q-1}(y_0,y_1)\geq 2^{q-1}$ and condition (\ref{eqn:Lem3})
holds.
\end{proof}

We now finish the proof of the theorem by handling the first
round, applying Lemma \ref{lem:manyPennies} until $x_0\leq 1$, and
by applying Lemma \ref{lem:smallMSB} until $ch(x_0,x_1)=0$.

\begin{proof}[Proof of Theorem \ref{thm:1Lie}]
By Lemma \ref{lem:1LieLB}, we may assume that $n$ satisfies
(\ref{eqn:1LieCondition}). For even $n=2m$, Paul chooses
$\vec{a}=(m,0)$ for his first
question so that the next state is forced to be $(y_0,y_1)=(m,m)$.
By Lemma \ref{lem:conserve} and the hypothesis, $wt_{q-1}(m,m)\geq
2^{q-1}$, and so $ch(m,m)\geq q-1$. If $m=1$, we apply Lemma
\ref{lem:smallMSB} to have Paul ask $q-1$ more questions.
Otherwise, $m>1$, and $(m,m)$ satisfies the requirements of Lemma
\ref{lem:manyPennies}.  We apply it repeatedly until we reach a
state of the form $(1,u)$. The lemma assures us that this will
happen in $t$ steps, where $\lfloor \log_2(m) \rfloor \le t \le
\lceil \log_2(m) \rceil$.  At the conclusion, we will have
$ch(1,u)\ge q-1-t$. Then, applying Lemma \ref{lem:smallMSB}, Paul
can ask at least $q-1-t$ further questions. Therefore, altogether
he has asked $1+t+(q-1-t)=q$ questions.

For odd $n=2m+1$, Paul chooses $\vec{a}=(m+1,0)$ for his first
question. Carole can then choose $(y_0,y_1)=(m+1,m)$ or $(m,m+1)$
as the next state.  We see that $wt_{q-1}(y_0,y_1)\ge
mq+m+1=\frac{2mq+2m+2}{2}=\frac{n(q+1)-(q-1)}{2}\ge 2^{q-1}$, by
hypothesis.  Hence regardless of Carole's response,
$ch(y_0,y_1)\ge q-1$.  The rest of the proof mimics the case for
even $n$.
\end{proof}

\section{Exact result for the 2-lie game\label{sec:twoLie}}

We now consider the $q$-round pathological liar game with 2 lies
and initial state $(n,0,0)$. The next theorem completely
characterizes the values of $n$ for which Paul can win the
$((n,0,0),q,2)^*$-game.   Its proof follows some definitions and
two lemmas focusing on the first two rounds and then the rest of the
game.

\begin{thm}\label{thm:2Lie}Let $q\geq 0$.  Paul has a winning strategy
for the $q$-round pathological liar game with 2 lies and initial
state $(n,0,0)$ iff
\begin{equation}
2^q \leq
    n\binom{q}{\leq 2} - A \binom{q-1}{2} - B \binom{q-2}{1},
\label{eqn:2LieCondition}
\end{equation}
where $A=n\mbox{ mod }2$ and
$$ B = \left\{\begin{array}{rl}
    0, & \mbox{if $n\equiv 0\mod 4$}, \\
    2\cdot (q\mbox{ mod }2), & \mbox{if $n\equiv 1\mod 4$}, \\
    (1-q^3)\mbox{ mod } 4, & \mbox{if $n\equiv 2\mod 4$}, \\
    (1+q^3)\mbox{ mod } 4, & \mbox{if $n\equiv 3\mod 4$}.
  \end{array}\right.
$$
\end{thm}

We say that Paul {\em survives} the first two rounds of the
$((n,0,0),q,k)^*$-game provided he has a strategy which guarantees
that the $(q-2)$-weight of the state after two rounds is at least
$2^{q-2}$ regardless of Carole's responses.
Let $\vec{a}$ be Paul's first question, and if Carole's
response is ``Y'' (``N''), then let Paul's second question be
$\vec{b^{\mathrm{Y}}}$ ($\vec{b^{\mathrm{N}}}$).  Then Paul can
survive the first two rounds iff
\begin{align}
2^{q-2}  \leq  \max_{\vec{a},\vec{b^{\mathrm{Y}}},
  \vec{b^{\mathrm{N}}}}
\min  \big\{ & wt_{q-2}(\mbox{Y}(\mbox{Y}((n,0,0),\vec{a}),
  \vec{b^{\mathrm{Y}}}),
  wt_{q-2}(\mbox{N}(\mbox{Y}((n,0,0),\vec{a}),
  \vec{b^{\mathrm{Y}}}), \nonumber \\
& wt_{q-2}(\mbox{Y}(\mbox{N}((n,0,0),\vec{a}),
  \vec{b^{\mathrm{N}}}),
  wt_{q-2}(\mbox{N}(\mbox{N}((n,0,0),\vec{a}),
  \vec{b^{\mathrm{N}}}) \big\}, \nonumber
\end{align}
where $\vec{a}$, $\vec{b^{\mathrm{Y}}}$, and $\vec{b^{\mathrm{N}}}$
must be legal questions when they are asked.  Now define weight
imbalances
\begin{eqnarray}
\Delta_{q-1} & := & \Delta_{q-1}((n,0,0),\vec{a}), \nonumber \\
\Delta_{q-2}^{\mathrm{Y}} & := & \Delta_{q-2}
    (\mbox{Y}((n,0,0),\vec{a}), \vec{b^{\mathrm{Y}}}), \quad\mbox{and}
    \nonumber \\
\Delta_{q-2}^{\mathrm{N}} & := & \Delta_{q-2}(\mbox{N}((n,0,0),
    \vec{a}), \vec{b^{\mathrm{N}}});  \nonumber
\end{eqnarray}
where without loss of generality, we choose the questions
$\vec{a}$, $\vec{b^{\mathrm{Y}}}$, and $\vec{b^{\mathrm{N}}}$
so that $\Delta_{q-1}$, $\Delta_{q-2}^{\mathrm{Y}}$, and
$\Delta_{q-2}^{\mathrm{N}}$ are {\em nonnegative} (for instance, by
replacing $\vec{a}$ with $\vec{x}-\vec{a}$).  By Lemma
\ref{lem:conserve} and (\ref{eqn:imbalance}), Paul can survive the
first two rounds of the $((n,0,0),q,2)^*$-game iff
\begin{eqnarray}
2^{q}  & \leq & wt_q(n,0,0) + \Delta,\qquad \mbox{where}\nonumber \\
\Delta & := & \max_{\vec{a},\vec{b^{\mathrm{Y}}},
  \vec{b^{\mathrm{N}}}}
\min  \big\{  \Delta_{q-1}+2\Delta_{q-2}^{\mathrm{Y}},
   \Delta_{q-1}-2\Delta_{q-2}^{\mathrm{Y}},
   \nonumber \\
& & \qquad\qquad
 -\Delta_{q-1}+2\Delta_{q-2}^{\mathrm{N}},
   -\Delta_{q-1}-2\Delta_{q-2}^{\mathrm{N}}\big\}. \label{eqn:delta}
\end{eqnarray}
We have reduced the problem to finding the value of $\Delta$
because, given a fixed first question $\vec{a}$, we may refer to Section
5 of \cite{G90} to compute $\vec{b^{\mathrm{Y}}}$ and $\vec{b^{\mathrm{N}}}$
minimizing $\Delta_{q-2}^{\mathrm{Y}}$ and $\Delta_{q-2}^{\mathrm{N}}$,
respectively.

\begin{lem}\label{lem:first2Rounds}
Let $q\geq 19$, $n\geq 2^q/\binom{q}{\leq 2}$, and let $\Delta$ be
defined as in (\ref{eqn:delta}) for the $q$-round pathological
liar game with 2 lies and initial state $(n,0,0)$. Then
$$\Delta= -A\binom{q-1}{2} -B \binom{q-2}{1},$$
where $A$ and $B$ are defined as in Theorem \ref{thm:2Lie}.
Furthermore, Paul's strategy achieving $\Delta$ guarantees at least
$(q-2)^2+\binom{q-2}{\leq 2}$ pennies after the first two rounds.
\end{lem}
\begin{proof}
Write $n=4p+r$ and $q-2=4l+s$, where $0\leq r,s<4$. We consider
cases of the initial state $(n,0,0)$ based on the values of $r$
and $s$. In each case, there is only one choice of $\vec{a}$
achieving $\Delta$ because any other choice of $\vec{a}$ results
in $-\Delta_{q-1}-2\Delta_{q-2}^{\mathrm{N}} <\Delta$ (recall
that, without loss of generality, $\vec{a}$,
$\vec{b^{\mathrm{Y}}}$ and $\vec{b^{\mathrm{N}}}$ are chosen to
make $\Delta_{q-1}$, $\Delta_{q-2}^{\mathrm{Y}}$ and
$\Delta_{q-2}^{\mathrm{N}}$ nonnegative). We give Paul's strategy
for achieving $\Delta$ by listing the questions $\vec{a}$,
$\vec{b^{\mathrm{Y}}}$ and $\vec{b^{\mathrm{N}}}$ in each case
explicitly. Guzicki proved that the choices below of
$\vec{b^{\mathrm{Y}}}$ and $\vec{b^{\mathrm{N}}}$ minimize
$\Delta_{q-2}^{\mathrm{Y}}$ and $\Delta_{q-2}^{\mathrm{N}}$; we
omit the details and refer the interested reader in Section 5 of
\cite{G90}. The calculations for the minimum $q$ for which all
questions are legal and for which the resulting states have at
least $(q-2)^2$ pennies are tedious but straightforward, and thus
omitted.

{\em Case $n=4p$}. Set  $\vec{a}=(2p,0,0)$ and
$\vec{b^{\mathrm{Y}}}=\vec{b^{\mathrm{N}}}=(p,p,0)$ to achieve
$\Delta=0$ with unique possible resulting state $(p,2p,p)$.  The
resulting state has $p\geq (q-2)^2+\binom{q-2}{\leq 2}$ pennies
when $q\geq 19$.

{\em Case $n=4p+1$}. Set $\vec{a}=(2p+1,0,0)$ and
$\vec{b^{\mathrm{Y}}}=(p+1,p,0)$ in each subcase, so that
$\Delta_{q-1}=\binom{q-1}{2}$ and two
possible resulting states are $(p+1,2p,p)$ and $(p,2p+1,p)$.
{\em Subcase $2\not|(q-2)$}.
Set $\vec{b^{\mathrm{N}}}=(p,p+1,0)$ to achieve
$\Delta=-\binom{q-1}{2}-2\binom{q-2}{1}$ with additional
possible resulting state $(p,2p,p+1)$.
{\em Subcase $2|(q-2)$}.  Set
$\vec{b^{\mathrm{N}}}=(p+1,p-\frac{q-2}{2}+1,0)$ to achieve
$\Delta=-\binom{q-1}{2}$ with additional possible resulting states
$(p+1,2p-\frac{q-2}{2},p+\frac{q-2}{2})$ and
$(p-1,2p+\frac{q-2}{2}+1,p-\frac{q-2}{2}+1)$.  All questions are legal
when $q\geq 9$, and all resulting states have at least
$(q-2)^2+\binom{q-2}{\leq 2}$ pennies when $q\geq 19$.

{\em Case $n=4p+2$}. Set $\vec{a}=(2p+1,0,0)$ in each subcase.
{\em Subcase $q-2=4l$}. Set
$\vec{b^{\mathrm{Y}}}=\vec{b^{\mathrm{N}}}=(p+1,p-\frac{q-2}{4}+1,0)$
to achieve $\Delta=-\binom{q-2}{1}$ with possible resulting states
$(p+1,2p-\frac{q-2}{4}+1,p+\frac{q-2}{4})$ and
$(p,2p+\frac{q-2}{4}+1,p-\frac{q-2}{4}+1)$.
{\em Subcase $q-2=4l+1$}. Set
$\vec{b^{\mathrm{Y}}}=\vec{b^{\mathrm{N}}}=(p,p+\frac{q-3}{4}+1,0)$
to achieve $\Delta=-2\binom{q-2}{1}$ with possible resulting states
$(p,2p+\frac{q-3}{4}+2,p-\frac{q-3}{4})$ and
$(p+1,2p-\frac{q-3}{4},p+\frac{q-3}{4}+1)$.
{\em Subcase $q-2=4l+2$}. Set
$\vec{b^{\mathrm{Y}}}=\vec{b^{\mathrm{N}}}=(p,p+\frac{q-4}{4}+1,0)$
to achieve $\Delta=-\binom{q-2}{1}$ with possible resulting states
$(p,2p+\frac{q-4}{4}+2,p-\frac{q-4}{4})$ and
$(p+1,2p-\frac{q-4}{4},p+\frac{q-4}{4}+1)$.
{\em Subcase $q-2=4l+3$}. Set
$\vec{b^{\mathrm{Y}}}=\vec{b^{\mathrm{N}}}=(p,p+\frac{q-5}{4}+1,0)$
to achieve $\Delta=0$ with possible resulting states
$(p,2p+\frac{q-5}{4}+2,p-\frac{q-5}{4})$ and
$(p+1,2p-\frac{q-5}{4},p+\frac{q-5}{4}+1)$. All questions are legal
when $q\geq 8$, and all resulting states have at least
$(q-2)^2+\binom{q-2}{\leq 2}$ pennies when $q\geq 19$.

{\em Case $n=4p+3$}. Set $\vec{a}=(2p+2,0,0)$ and
$\vec{b^{\mathrm{Y}}}=(p+1,p+1,0)$ in each subcase so that two
possible resulting states are always $(p+1,2p+2,p)$ and
$(p+1,2p+1,p+1)$.
{\em Subcase $q-2=4l$}. Set
$\vec{b^{\mathrm{N}}}=(p,p+\frac{q-2}{4}+1,0)$
to achieve $\Delta=-\binom{q-1}{2}-\binom{q-2}{1}$
with additional possible resulting states
$(p,2p+\frac{q-2}{4}+2,p-\frac{q-2}{4}+1)$ and
$(p+1,2p-\frac{q-2}{4}+1,p+\frac{q-2}{4}+1)$.
{\em Subcase $q-2=4l+1$}. Set
$\vec{b^{\mathrm{N}}}=(p,p+\frac{q-3}{4}+1,0)$
to achieve $\Delta=-\binom{q-1}{2}$
with additional possible resulting states
$(p,2p+\frac{q-3}{4}+2,p-\frac{q-3}{4}+1)$ and
$(p+1,2p-\frac{q-3}{4}+1,p+\frac{q-3}{4}+1)$.
{\em Subcase $q-2=4l+2$}. Set
$\vec{b^{\mathrm{N}}}=(p+1,p-\frac{q-4}{4}+1,0)$
to achieve $\Delta=-\binom{q-1}{2}-\binom{q-2}{1}$
with additional possible resulting states
$(p+1,2p-\frac{q-4}{4}+1,p+\frac{q-4}{4}+1)$ and
$(p,2p+\frac{q-4}{4}+2,p-\frac{q-4}{4}+1)$.
{\em Subcase $q-2=4l+3$}. Set
$\vec{b^{\mathrm{N}}}=(p,p+\frac{q-5}{4}+2,0)$
to achieve $\Delta=-\binom{q-1}{2}-2\binom{q-2}{1}$
with additional possible resulting states
$(p,2p+\frac{q-5}{4}+3,p-\frac{q-5}{4})$ and
$(p+1,2p-\frac{q-5}{4},p+\frac{q-5}{4}+2)$.
All questions are legal
when $q\geq 8$, and all resulting states have at least
$(q-2)^2+\binom{q-2}{\leq 2}$ pennies when $q\geq 19$.
\end{proof}

\begin{lem}\label{lem:2LiePennies}
Let $q \geq 23$. If $wt_q(x_0, x_1, x_2)=2^q$ and $x_2 \geq q^2$,
then Paul has a strategy to reach a state $\vec z$ with
$wt_0(\vec{z})=1$ in exactly $q$ rounds such that every
intermediate state $(u_0, u_1, u_2)$ after playing $q-j$ rounds
satisfies $wt_j(u_0,u_1,u_2)=2^j$.
\end{lem}
\begin{proof}
The proof proceeds by showing how Spencer's ``Main Theorem'' of
\cite[Section 2]{S92}, quoted here as Theorem \ref{thm:spencer},
can be tightened in the case $k=2$ so that we may take $c=1$ and
$q_0=23$. Spencer's technique is to relax the game to allow the
pennies position to take on negative integer values in both
questions and resulting states in {\em fictitious play}, and then
to show in fact that this position never goes negative for a given
$c$ and $q_0$.

Before stating and proving the three claims which tighten Spencer's
result, we recall the necessary notation and results from \cite{S92}
for the case $k=2$.  Assume there are $j+1$ rounds remaining,
and the current position is $\vec{P}=(p_0,p_1,p_2)$ with weight $2^{j+1}$.

{\em Fictitious play}:  Paul  selects the next question vector
$(v_0,v_1,v_2)$ according to the parity of $p_0$ and $p_1$ as
follows. If $p_0$ is odd, then $v_0=\frac{p_0+1}{2}$ and
$v_1=\lfloor\frac{v_1}{2}\rfloor$; otherwise if $p_0$ is even,
then $v_0=\frac{p_0}{2}$ and $v_1=\lceil\frac{v_1}{2}\rceil$.
Let $v_2$ be the unique integer that makes the weight imbalance
$\Delta_j(\vec P, \vec v)=0$.  In other words, in fictitious
play the weight of the states is exactly halved after each round.
Note that by the choices of $v_0$ and $v_1$,
$\Delta_j(\vec P, \vec v)-(2v_2-p_2)\geq 0$. Hence,  $v_2 \leq p_2$,
and so $(v_0,v_1,v_2)$ is legal whenever $v_2\geq 0$.

In fictitious play Paul and Carole continue to play formally even
though the last entry of the states may turn negative. Let
\[
fic(j) =(fic_0(j), fic_1(j), fic_2(j))
\]
be the state of the game when there are $j$ rounds remaining. Note
that $fic(q)=(x_0, x_1, x_2)$ is simply the initial state of the
game, and $fic_0(j), fic_1(j)$ are always non-negative.

{\em Perfect play}: When the state is $\vec P$, Paul selects
$\vec v=\vec P /2$. This results in $Y(\vec P, \vec v)=N(\vec P,
\vec v)$ and uniquely determines the state $pp(j)=(pp_0(j),
pp_1(j), pp_2(j))$ when $j$ rounds remain in the game. When
the initial state $pp(q)$ is  $\vec{x}$, it is easy
to compute that
\begin{eqnarray*}
pp_0(j)& =& \frac{x_0}{2^{q-j}}, \qquad pp_1(j)=\frac{ x_1+x_0
{q-j \choose 1}}
{2^{q-j}}. \\
pp_2(j) & = & \frac{x_2+ x_1 {q-j \choose 1} +x_2 { q-j \choose
2}}{2^{q-j}}.
\end{eqnarray*}

Defining $e_i(j)=|pp_i(j)-fic_i(j)|$,
Spencer proves $e_0(j)\leq 1$ and $e_1(j) \leq 3$.  By replacing
the $j^k$ in Spencer's calculations with $\frac{1}{2}\binom{j}{k}$
for $k=2$, it follows that
$|fic_2(j)-\frac{1}{2}(fic_2(j+1)+fic_1(j+1))| \leq \frac{1}{2} {j
\choose 2} +\frac{1}{2}$. Hence $e_2(j)\leq \frac{1}{2}e_2(j+1)
+\frac{1}{2}{j \choose 2} +2$. By induction, $e_2(j) \leq {j
\choose 2}+5$.

Now we describe the strategy for Paul: starting from the state
$(x_0, x_1, x_2)$ with $q$-weight $2^q$ and $x_2 \geq q^2$,
Paul plays fictitious play in all rounds.  Our analysis now
deviates from that of Spencer.
We argue that Paul can win by seeing that no entries turn negative
and by examining the state Paul reaches at
$j=6$, i.e., when 6 rounds remain. Explicitly, we prove
the following claims for $q \geq 23$.\smallskip

\begin{enumerate}
\item[1.] $fic_0(6) \leq 1$.
\item[2.]  $fic_2(j) >1$ for $j \geq 6$. (Fictitious
play questions are legal when $j \geq 6$.)
\item[3.] When $j=6$, the state of the game is not $(1, 5,
7)$, or $(1, 4, 14)$.
\end{enumerate}

If the above claims are true, then the possible states at $j=6$
are $(1,3,21)$, $(1,2,18)$, $(1,1,35)$, $(1,0,42)$, $(0,8,8)$,
$(0,7,15)$, $(0,6,22)$, $(0,5,29)$, $(0,4,36)$, $(0,3,43)$,
$(0,2,50)$, $(0,1,57)$, and $(0,0,64)$. It is easy to check that in
all these states, Paul can split the weight evenly until he
reaches a state $\vec z$ with $wt_0(\vec z)=1$.

{\bf Proof of Claim 1}. \ Since $e_0(j)\leq1$, it suffices to show that
$pp_0(6) <1$, i.e., $x_0 < 2^{q-6}$. This is true  because  $x_0
{q \choose \leq 2} \leq wt_q(x_0,x_1,x_2) =2^q$. Hence $x_0 \leq
2^q /{q \choose \leq 2}$, which is less than $2^{q-6}$ when $q
\geq 12$.

{\bf Proof of Claim 2}. \ We show that $pp_2(j) >e_2(j)+1$ for
$6\leq j\leq q-1$. It is
enough to show that $\min\left\{\left(x_2+x_1\binom{q-j}{1}
+x_0{q-j\choose 2}\right)/2^{q-j}\right\}  \geq {j \choose 2}+6 $
for all $x_0, x_1, x_2$ satisfying $wt_q(x_0,x_1,x_2)=2^q$
and $x_2 \geq q^2$. The minimum of $\left(x_2+x_1(q-j)
+x_0{q-j\choose 2}\right)/2^{q-j}$ is achieved at one of the
vertices of the feasible region, that is, when $(x_0, x_1, x_2)$
is
\[
(0,0,2^q), \ \left(0, \frac{2^q-q^2}{\binom{q}{\leq 1}},
q^2\right), \text{ or }
\left(\frac{2^q-q^2}{{q \choose \leq 2}},0,q^2\right).
\]
For $6 \leq j \leq q-1$, direct computation shows that the minimum
is greater than ${j \choose 2}+6 > e_2(j)+1$ for $q \geq 16$. The
case $j=q-1$ is special, for which ${q-j \choose 2}=0$;
the inequality remains true here since $x_2 \geq q^2$.

{\bf Proof of Claim 3.} \ We show that when $q\geq 23$,
$fic_0(8)\leq 1$ and $fic_0(8)+fic_1(8) \leq 16$. Then by
definition of fictitious play, $fic(6)$ could not be $(1, 5, 7)$ or
$(1, 4, 14)$.

To show $fic_0(8) \leq 1$, note that $pp_0(8)=x_0/2^{q-8} \leq
2^8/{q\choose 2}$, which is less than 1 when $q \geq 23$. To show
$fic_0(8)+fic_1(8) \leq 16$, define
$e_{01}(j):=|pp_0(j)+pp_1(j)-(fic_0(j)+fic_1(j))|$. By definition
of fictitious play,
$|fic_0(j)+fic_1(j)-(fic_0(j+1)+\frac{1}{2}fic_1(j+1))| \leq
\frac{1}{2}$, which implies $e_{01}(j)\leq \frac{1}{2}
e_{01}(j+1)+1$. By induction with base case $e_{01}(q)=0$, we have
$e_{01}(j) < 2$. Now assume that $fic_0(8)+fic_1(8) \geq 17$. Then
$pp_0(8)+pp_1(8)
>15$, that is, $x_1+x_0\binom{q-8}{\leq 1} > 15\cdot 2^{q-8}$.
However, the maximum of $x_1+x_0\binom{q-8}{\leq 1}$ is reached
when $(x_0, x_1)$ is either $(0, (2^q-q^2)/\binom{q}{\leq 1})$, or
$((2^q-q^2)/{q\choose \leq 2}, 0)$. For the first one,
$x_1+x_0\binom{q-8}{\leq 1}> 15\cdot 2^{q-8}$ iff $q \leq 16$, for
the second one, $x_1+x_0\binom{q-8}{\leq 1}
> 15\cdot 2^{q-8}$ iff $q \leq 22$. This contradicts the fact
that $q \geq 23$.
\end{proof}

\begin{proof}[Proof of Theorem \ref{thm:2Lie}]
The values of $F^*_2(q)$ for $1\leq q\leq 24$, found by exhaustive
computation, are listed in Table \ref{tab:2Lie}.  In each case,
$F^*_2(q)$ is the first value of $n$ which satisfies the
inequality in (\ref{eqn:2LieCondition}). These values were
generated by a dynamic programming algorithm based on the
recurrence
$$ r^*(\vec{x}) = 1 + \max_{\vec{a}}\{\min\{r^*(\mbox{Y}(\vec{x},\vec{a})),
    r^*(\mbox{N}(\vec{x},\vec{a}))\}\},  $$
where $r^*(\vec{x})$ is defined to be the maximum number of rounds
for which Paul can win the pathological liar game with initial
state $\vec{x}$.

\begin{table}[h]
\begin{center}{\small
\begin{tabular}{|r||c|c|c|c|c|c|c|c|c|c|c|c|}
\hline $q$ & 1 & 2 & 3 & 4 & 5 & 6 & 7 & 8 & 9 & 10 & 11 & 12\\
\hline $F^*_2(q)$ & 1 & 1 & 2 & 2 & 2 & 4 & 6 & 8 & 12 & 20 & 32 & 52\\
\hline
\multicolumn{13}{c}{\phantom{m}}\\
 \hline $q$ & 13 & 14 & 15 & 16 & 17 & 18 & 19 & 20 & 21 &
22 & 23 & 24 \\
\hline $F^*_2(q)$ & 90 & 156 & 272 & 480 & 852 & 1525 & 2746 &
4970 & 9040 & 16514 & 30284 & 55740 \\ \hline
\end{tabular}}
\end{center}
\caption{Values of $F^*_2(q)$, the minimum number of elements n
for which Paul can win the $q$-round pathological liar game with 2
lies and initial state $(n,0,0)$. \label{tab:2Lie}}
\end{table}

Now suppose $q\geq 25$.  If $n$ satisfies
(\ref{eqn:2LieCondition}), then by Lemma \ref{lem:first2Rounds},
Paul can survive two rounds with all possible resulting states
having at least $(q-2)^2+\binom{q-2}{\leq 2}$ pennies. If after
the first two rounds the $(q-2)$-weight of the resulting state is
$>2^{q-2}$, greedily remove coins as large as possible so that
the $(q-2)$-weight is exactly $2^{q-2}$.  Either the resulting state
has only pennies
remaining, or at most $\binom{q-2}{\leq 2}$ pennies were removed.
Since $q-2\geq 23$, Lemma
\ref{lem:2LiePennies} shows that Paul can win the
$((n,0,0),q,2)^*$-game. If $n$ fails to satisfy
(\ref{eqn:2LieCondition}), then by (\ref{eqn:delta}) and Lemma
\ref{lem:first2Rounds}, Paul cannot survive the first two rounds
and therefore has no winning strategy for the
$((n,0,0),q,2)^*$-game.
\end{proof}

\section{Winning strategies and hypercube coverings and packings
\label{sec:equiv}}

The pathological liar game has an important natural reformulation
in terms of coverings of the hypercube $Q_k$ with certain adaptive
Hamming balls. For our purposes, we think of the $q$-dimensional
hypercube $Q_q$ as the set of vertices $\{\mbox{Y},\mbox{N}\}^q$ in
which two vertices are adjacent iff they differ in exactly one
position. Instead of the usual 0's and 1's, the bits are Y's and
N's, and so ``bit'' complementation is defined by
$\overline{\mbox{Y}}=\mbox{N}$ and $\overline{\mbox{N}}=\mbox{Y}$.
A Hamming ball of radius $k$ in $Q_q$ consists of a center
$\omega\in Q_q$ and all $w'\in Q_q$ which differ from $\omega$ in
at most $k$ positions. A {\em covering} ({\em packing}) of $Q_q$
usually refers to a collection of Hamming balls of a fixed radius
whose union is $Q_q$ (disjoint in $Q_q$), but there are many
variations.  We refer the interested reader to the literature for
further information \cite{CHLL97,BHP98}.
It happens that a winning strategy for Paul in the pathological
liar game can be converted to a covering of $Q_q$ with these
adaptive Hamming balls, and vice versa.  We now formalize this
relationship.

\begin{figure}[h]
\begin{center}
\psfig{figure=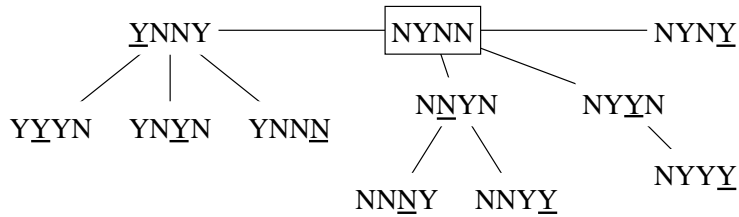}
\caption{A 2-quasiball in the hypercube $Q_4$ satisfying Definition
\ref{def:iball} is given a tree structure with stem NYNN. A child
agrees with its parent before the underlined position, is opposite
at the underlined position, and has unconstrained relationship
with its parent afterward.
\label{fig:iball}}
\end{center}
\end{figure}

Noting that $2^{[q]}$ is the power set of $[q]$, define
$$\binom{[q]}{j}:=\{J\in 2^{[q]}:|J|=j\} \qquad \mbox{and} \qquad
\binom{[q]}{\leq i} \ := \ \bigcup_{j=0}^i \binom{[q]}{j}. $$
We have the following definition of an adaptive Hamming ball, which
we call a {\em quasiball}, followed by an example for $q=4$ and
radius $i=2$.
\begin{defn}[$i$-quasiball]\label{def:iball}
Let $q,i\geq 0$. An $i$-quasiball is the image
$f\left(\binom{[q]}{\leq i}\right)$ of an injective function
$$
f \ : \ \binom{[q]}{\leq i} \rightarrow Q_q,  $$
such that whenever $A,B\in \binom{[q]}{\leq i}$ are of the form
\begin{equation}\label{eqn:ABForm}
A = \{p_1,\ldots,p_{|A|}\} \quad \mbox{and} \quad
   B = \{p_1,\ldots,p_{|A|},p_{|A|+1},\ldots,p_{|B|}\},
\end{equation}
where $p_1< \cdots <p_{|A|}<p_{|A|+1}<\cdots<p_{|B|}$, then
$f(A)$ and $f(B)$ are of the form
\begin{equation}\nonumber
f(A) = \omega_1\cdots \omega_{p_{|A|}}\cdots \omega_q \quad \mbox{and}
    \quad f(B) = \omega_1\cdots \omega_{(p_{|A|+1}-1)}
    \overline{\omega}_{p_{|A|+1}}
    \omega'_{(p_{|A|+1}+1)}\cdots\omega'_q,
\end{equation}
where $\omega'_{(p_{|A|+1}+1)}\cdots\omega'_q \in Q_{q-p_{|A|+1}}$.
\end{defn}
\begin{exmp}[A 2-quasiball in $Q_4$] Let $q=4$ and $i=2$.
Define $f:\binom{[4]}{\leq 2}\rightarrow Q_4$ by
$f(\emptyset)=\mbox{NYNN}$, $f(\{1\})=\mbox{YNNY}$,
$f(\{2\})=\mbox{NNYN}$, $f(\{3\})=\mbox{NYYN}$,
$f(\{4\})=\mbox{NYNY}$, $f(\{1,2\})=\mbox{YYYN}$,
$f(\{1,3\})=\mbox{YNYN}$, $f(\{1,4\})=\mbox{YNNN}$,
$f(\{2,3\})=\mbox{NNNY}$, $f(\{2,4\})=\mbox{NNYY}$, and
$f(\{3,4\})=\mbox{NYYY}$.  For instance, letting $A=\{2\}$
and $B=\{2,3\}$, we see that the first two coordinates of $f(A)$
and $f(B)$ agree, and the third coordinate is opposite, satisfying
the constraint on $A$ and $B$ given by the definition (the fourth
coordinate happens to be opposite as well).  After similar verification
for all possible choices of $A$ and $B$, we see that
 $f\left(\binom{[4]}{\leq 2}\right)$ is a $2$-quasiball in
$Q_4$. We assign a tree structure to
$f\left(\binom{[4]}{\leq 2}\right)$ by defining the parent of $f(B)$,
for any $B=\{p_1,\ldots,p_{|B|}\}\neq \emptyset$,
to be $f(B\setminus\{p_{|B|}\})$, as illustrated in Figure
\ref{fig:iball}.
\end{exmp}
Intuitively, for $i>0$ an $i$-quasiball contains a stem, $f(\emptyset)$,
and $q$ children in $f\left(\binom{[q]}{1}\right)$ obtained from
$f(\emptyset)$  by complementing one of its $q$ bits and
choosing  the bits to the right arbitrarily.  The child $f(\{p\})$
can be considered to be the stem of the $(i-1)$-quasiball obtained by
deleting the first $p$ bits from each of the vertices in
$f(\{\{p\}\cup P : P\in \binom{[q]\setminus [p]}{\leq i-1}\})$. An
$i$-quasiball is clearly a generalization of a Hamming ball of radius
$i$, since for $A,B\in Q_q$ satisfying (\ref{eqn:ABForm}), we may
choose $f(B)$ by complementing $f(A)$ in positions
$p_{|A|+1},\ldots,p_{|B|}$ and leaving the other positions
unchanged. We note in passing that some $i$-quasiballs, for example
$\{\mbox{Y},\mbox{N}\}$ and otherwise whenever $i\geq q$, are
obtained from more than one such function $f$.

In order to understand the relationship between winning strategies
for Paul and coverings by $i$-quasiballs, recall that a {\em covering
code} of length $q$ and radius $k$ is a set of Hamming balls of
radius $k$ whose union is $Q_q$.  By relaxing Hamming balls to
$i$-quasiballs and by allowing $i$ to vary between $0$ and $k$, we
define an {\em $\vec{x}$-covering}, where
$\vec{x}=(x_0,\ldots,x_k)$ to be a collection consisting of $x_i$
$(k-i)$-quasiballs for each $0\leq i\leq k$ whose union is $Q_q$.
Similarly, an {\em $\vec{x}$-packing} is such a collection whose
constituent members are pairwise disjoint, and whose union is not
necessarily $Q_q$. Informally speaking, we may think of an
$(n,0,\ldots,0)$-covering of $Q_q$ as an adaptive covering code
of length $q$ and fixed radius.
The following theorem is adapted from
\cite[Theorem 1.2]{SY03} which is for an asymmetric version of the
original game.
\begin{thm}\label{thm:equivCovering}
Let $q,k\geq 0$.  Paul has a strategy for winning the $q$-round,
$k$-lie pathological liar game with initial state $\vec{x}$ iff
there exists an $\vec{x}$-covering of $Q_q$.  Similarly, Paul has
a strategy for winning the original game with the same parameters
iff there exists an $\vec{x}$-packing of $Q_q$.
\end{thm}
\begin{proof}
For the proof it is convenient to keep track of the sets of
elements with a given number of lies, and not just their
cardinalities. Without loss of generality, in a game with initial
state $\vec{x}=(x_0,\ldots,x_k)$, let $n=\sum_{i=0}^k x_i$ and let
$X_i\subseteq [n]$ be the $x_i$ elements initially associated with
$i$ lies. We will abuse notation and let a state or question
vector be given in either integer or set format; for example,
$\vec{x}=(x_0,\ldots,x_k)$ or $(X_0,\ldots,X_k)$. We prove
the statement about the pathological liar game and remark how to
adapt the proof for the original game afterward.

For the forward implication, Paul's winning strategy corresponds
to a decision tree which is a full binary tree of depth $q$.  The
root contains the initial state $\vec{x}$ and the first question.
Each node contains a nonzero state, and each internal node
contains a legal question for the state in the same node.  A node
containing state $\vec{P}$ and question $\vec{v}$ has left child
containing state $\mbox{N}(\vec{P},\vec{v})$ and right child
containing state $\mbox{Y}(\vec{P},\vec{v})$, corresponding to
responses of ``N'' or ``Y,'' respectively, by Carole. A game
played under this strategy is a path from the root to a leaf of
the decision tree, passing down $q$ levels of questions by Paul
and answers by Carole. We say that a leaf is {\em labeled} by each
element of $[n]$ which survives in that leaf's state. A leaf
labeled by $x\in[n]$ has a {\em response vertex} with respect to
$x$, which is Carole's Yes/No response sequence
$\omega_1\cdots\omega_q\in Q_q$ read in order from the root to
that leaf.  If the context is clear, we will refer to a response
vertex with respect to $x$ simply as a response vertex. The leaves
are in bijection with $Q_q$ by considering the response sequence
leading to each leaf.

Let $i\in\{0,\ldots,k\}$ and choose $x\in X_i$. Let $S\subseteq
Q_q$ be the set of response vertices with respect to $x$ of those
leaves labeled with $x$. We define the function
$f:\binom{[q]}{\leq i}\rightarrow Q_q$ certifying that $S$ is a
$(k-i)$-quasiball as follows.  Set $f(\emptyset)$ equal to the unique
$\omega\in S$ for which every response by Carole is truthful.  In
general a response vertex is completely determined by the
positions $A\subseteq [q]$ corresponding to lies by Carole.  Set
$f(A)$ equal to this response vertex for all $A\in
\binom{[q]}{\leq k-i}$. Two leaves $\alpha$ and $\beta$ both
labeled by $x$ and having response vertices with lies in positions
$A,B\subseteq [q]$, respectively, and satisfying
(\ref{eqn:ABForm}), must have the same first $p_{|A|}-1$ response
sequence steps from the root and bifurcate at step $p_{|A|}$.
Therefore $S$ is a $(k-i)$-quasiball, and since every leaf is labeled
by at least one element of $[n]$, there exists an
$\vec{x}$-covering of $Q_q$.

For the reverse implication, the states and questions contained in
the depth $q$ full binary decision tree are determined by the
$\vec{x}$-covering. The initial state at the root is
$\vec{x}=(X_0,\ldots,X_k)$, where each $(k-i)$-quasiball is identified
with a unique element $x\in X_i$ and is the image of a function
$f_x:\left(\binom{[q]}{\leq k-i}\right)\rightarrow Q_q$ satisfying
Definition \ref{def:iball}. Paul constructs the first question
vector $\vec{a}=(A_0,\ldots,A_k)$ by letting $x\in A_i$ whenever
the stem of the $(k-i)$-quasiball identified with $x$ begins with
``Y.''  Thus every $x\in A_i$ will label a leaf whose response
vertex with respect to $x$ begins with ``Y.'' Suppose Carole
responds to $\vec{a}$ with ``Y.'' If $x\in A_i$ for some $i$, no
lie is associated with $x$ by Carole's response, and
$f_x\left(\binom{[q]\setminus \{1\}}{\leq k-i}\right)\subseteq
\mbox{Y}Q_{q-1}$ may be viewed as a $(k-i)$-quasiball in $Q_{q-1}$
by restricting the domain of $f_x$ to $\binom{[q]\setminus
\{1\}}{\leq k-i}$ and deleting the first bit of each vertex in the
image.  The resulting state vector $\mbox{Y}(\vec{x},\vec{a})$
counts $x$ in the $i$th position.  If $x$ is not counted by
$\vec{a}$, one lie is associated to $x$ by Carole's response, and
$f_x(\{\{1\}\cup P:P\in \binom{[q]\setminus \{1\}}{\leq
k-i-1}\})\subseteq \mbox{Y}Q_{q-1}$ may be viewed as a
$(k-i-1)$-quasiball in $Q_{q-1}$ by restricting the domain of $f_x$ to
$\{\{1\}\cup P:P\in \binom{[q]\setminus \{1\}}{\leq k-i-1}\}$ and
deleting the first bit of each vertex in the image.  The resulting
state vector  $\mbox{Y}(\vec{x},\vec{a})$ counts $x$ in the
$(i+1)$st position (if $i+1>k$, then the $(k-i-1)$-quasiball is empty
and $x$ does not appear in $\mbox{Y}(\vec{x},\vec{a})$).
In both cases, the rest of the domain of $f_x$ is mapped to
$\mbox{N}Q_{q-1}$. Therefore
there exists a $\mbox{Y}(\vec{x},\vec{a})$-covering of $Q_{q-1}$.
Similarly, if Carole answers ``N'' there exists a
$\mbox{N}(\vec{x},\vec{a})$-covering of $Q_{q-1}$.  The reverse
implication follows by induction, since a covering of $Q_0$ must
consist of at least one $i$-quasiball, which corresponds to a surviving
element.

For the original liar game, the function $f$ in the forward
implication is defined in the same way; however, there is at most
one surviving element labeling each leaf of the decision tree.
This ensures that the collection of $i$-quasiballs which are the sets
of response vertices of leaves with a given label are disjoint,
and thus form a packing.  For the reverse implication, the
inductive step is the same, but for the base case a packing of
$Q_0$ corresponds to at most one $i$-quasiball.
\end{proof}

Monotonicity under majorization, defined in Section \ref{sec:vector},
is now clear because an $i$-quasiball realized by a function
$f:\binom{[q]}{\leq i}\rightarrow Q_q$ can be considered to contain
an $(i-1)$-quasiball obtained by restricting $f$ to
$\binom{[q]}{\leq i-1}$.
Theorem \ref{thm:equivCovering} allows
Lemma \ref{lem:lowerBound}, and its dual version
for the original game, to be interpreted in terms of the sphere
bound for coverings or packings, respectively, of the hypercube.
A $k$-quasiball has size $\binom{q}{\leq k}$ in $Q_q$, and
so there can exist neither a covering of $Q_q$ with fewer than
$2^q/\binom{q}{\leq k}$ $k$-quasiballs, nor a packing of $Q_q$ with
more than $2^q/\binom{q}{\leq k}$ $k$-quasiballs.
A natural question is whether
the asymptotic sizes of optimal coverings and packings, that is,
covering codes and error-correcting codes, meet at the
sphere bound.  For Hamming balls, this is true for radius 1
\cite[Theorem 12.4.11]{CHLL97}, and is
unknown for larger radius. For $k$-quasiballs,
this is now known to be true for fixed $k$ by combining Theorems
\ref{thm:spencerMain} and \ref{thm:asymptotic}.

\section*{Acknowledgment}

We would like to thank Joel Spencer for a helpful discussion with
the first author.


\end{document}